\RequirePackage{fix-cm}
\documentclass[smallextended]{svjour3}       
\smartqed  
\usepackage{lineno,hyperref}
\usepackage{graphicx}
\usepackage{lineno,hyperref}
\usepackage{amsmath}
\usepackage{amsfonts}
\usepackage{amssymb}
\usepackage{graphicx}
\usepackage{caption}
\usepackage{subcaption}
\usepackage{float}
\usepackage{multirow}
\usepackage{boldline,multirow}
\usepackage{tabularx,ragged2e}
\usepackage{boldline}
\usepackage{tabu}
\usepackage{longtable}
\usepackage{color}
\begin{document}

\title{A comparison between pre-Newton and post-Newton approaches for solving a physical singular second-order boundary problem in the semi-infinite interval}

\titlerunning{Tow efficient numerical methods for Thomas-Fermi equation}        

\author{Amir Hosein Hadian-Rasanan*\and
      Mehran Nikarya \and
      Arman Bahramnezhad \and
      Mohammad M. Moayeri\and
      Kourosh Parand
}


\institute{A. H. Hadian-Rasanan \at
				Department of Computer Sciences, Shahid Beheshti University, G.C., Tehran, Iran.\\
				\email{amir.h.hadian@gmail.com}
         \and
         M. Nikarya \at
          Department of Computer Sciences, Shahid Beheshti University, G.C., Tehran, Iran.\\
          \email{mehran.nikarya@gmail.com}
         \and
         A. Bahramnezhad \at
            Department of Computer Sciences, Shahid Beheshti University, G.C., Tehran, Iran.\\
            \email{arman.bahramnezhad@gmail.com}
         \and
         M. M. Moayeri \at
         Department of Computer Sciences, Shahid Beheshti University, G.C., Tehran, Iran.\\
         \email{m\_moayeri@sbu.ac.ir}
         \and
         K. Parand \at
         Department of Computer Sciences, Shahid Beheshti University, G.C., Tehran, Iran. \\
         Department of Cognitive Modeling, Institute for Cognitive and Brain Sciences, Shahid Beheshti University, G.C., Tehran, Iran.\\
         \email{k\_parand@sbu.ac.ir}
}

\date{Received: date / Accepted: date}

\maketitle

\begin{abstract}
In this paper, two numerical approaches based on the Newton iteration method with spectral algorithms  are introduced to solve the Thomas-Fermi equation. That Thomas-Fermi equation is a nonlinear singular ordinary differential equation (ODE) with boundary condition in infinite. In these schemes, the Newton method is combined with a spectral method where in one of those, by Newton method we convert nonlinear ODE to a sequence of linear ODE then, solve them using the spectral method. In another one, by the spectral method the nonlinear ODE be converted to system of nonlinear algebraic equations, then, this system is solved by Newton method. In both approaches, the spectral method is based on fractional order of rational Gegenbauer functions. Finally, the obtained results of two introduced schemes are compared to each other in accuracy, runtime and iteration number. Numerical experiments are presented showing that our methods are as accurate as the best results which obtained until now.
\end{abstract}
\section{Introduction}
\paragraph{}
Nonlinear problems arise in the various field of research such as biology, cognitive sciences,
engineering, finance, etc. One of the branches of nonlinear problems is nonlinear ODEs which have
unbounded domains. Since these problems are significant, many researchers developed different
numerical schemes to solve them. There are various numerical algorithms to compute solution of
nonlinear problems over the semi-infinite domains such as
Adomian decomposition method\cite{Adomian1998,Wazwaz1999,Randolph,Epele1999}, finite difference and finite element
methods\cite{Noye1999,Bu2015,Choi2016}, Hermite collocation\cite{Bayatbabolghani2014}, meshless
methods\cite{Hemami,Parand2011,Kazem2012}, etc. In this work, two different approaches based on the
combination of spectral methods and Newton family algorithms are introduced. As the first one, we can refer to the pre-Newton. In the pre-Newton approach, a linearization method is done directly on the nonlinear ODE; then, nonlinear ODE is converted  to a sequence of linear ODE which can be solved by different numerical algorithms such as spectral methods. The second approach is the post-Newton. In the post-Newton approach, by using a numerical method the nonlinear ODE is converted to a system of nonlinear algebraic equations and then, this system is solved using various Newton type algorithms.
To show the efficiency of these approaches and compare them to each other, we
consider a nonlinear ODE called Thomas-Fermi equation which arises in theoretical physics as a test
problem. This model has two significant roles in mathematical physics for two reasons: Thomas-Fermi equation was
enhanced to model the effective nuclear charge in heavy atoms, and  was investigated to analyze the
potentials and charge densities of atoms having numerous electrons \cite{Zhu2012}.

In this paper, pre-Newton and post-Newton approaches based on the fractional order of rational Gegenbauer (FRG) functions are used to solve the Thomas-Fermi equation in the semi-infinite interval. The main aim of this paper is presenting a kind collocation method based on FRG for the solving Thomas-Fermi equation which can obtain the most accurate results which are reported until now. In this paper, we are going to compute $y'(0)$. The obtained value for $y'(0)$ is as follows:$$−1.588071022611375312718684509423950109452746621674825616765677.$$ This value is obtained in \cite{Zhang2018} using 600 basis functions. But in this paper, we obtain this value by using only 200 basis function The main advantage of the presented method is highly convergence rate of it. On the other hand, it has a good time efficiency. The organization of the paper is expressed as follows: Thomas-Fermi equation is introduced in Section 2, The Gegenbauer polynomials and FRG functions are introduced in Section 3. Section 4 contains the Newton-Kantorovich method and the application of spectral methods. Results and discussion of the proposed methods are shown in Section 4. Finally, a conclusion is provided in Section 5.

\section{Thomas-Fermi equation}
\paragraph{}
The Thomas-Fermi theorem illustrates that how the energy of an electronic system, $E$, and the electronic density, $\rho$, are connected to each other by the following formula\cite{Parand2013}:
\begin{equation}
E[\rho] = \frac{9}{10B}\int\rho(r)d\tau + \frac{1}{2}\int\frac{\rho(r)\rho(r')}{|r-r'|}d\tau'd\tau + \int\rho(r)\nu(r)d\tau,
\end{equation}
where $\nu(r)$ is the external potential and $B = 3(3\pi)^{-\frac{2}{3}}$. In order to obtain the density the energy functional should be minimize with respect to $\rho$ and subject to the normalization restriction $\int \rho(r)d\tau = N$ where $N$ is the number of electrons
\begin{equation}
\frac{3}{2B}\rho(r)^{\frac{3}{2}} + \int\frac{\rho(r')}{|r-r'|}d\tau' + \nu(r) = \mu,
\end{equation}
where $\mu$ is the Lagrange multiplier related to the normalization restriction\cite{Cedillo1993,Parand2013}. By using Poisson's equation to remove the density and a change of variables Thomas-Fermi equation is obtained as below:
\begin{equation}
\label{TFE}
\frac{d^2y}{dx^2}=\frac{1}{\sqrt{x}}y^{\frac{3}{2}}(x),
\end{equation}
with the following boundary conditions:
\begin{eqnarray}
y(0)=1,~~~~\lim_{x\to\infty} y(x)=0.
\end{eqnarray}
This equation describes the charge density in atoms of high atomic number and appears in the problem of determining the effect of nuclear charge in heavy atoms\cite{Thomas1927,Davis1962,Parand2013}.

As the solution of Thomas-Fermi equation is effective in theoretical physics, many scientists has studied this model. Moreover this equation has three different forms which can be effective on the rate of convergence of the using numerical algorithm\cite{Zhang2018}. These three forms of Thomas-Fermi equation are listed in Table \ref{formstable}.
\begin{table}
\centering
\caption{Different forms of Thomas-Fermi equation}
\begin{tabularx}{\textwidth}{p{4cm} p{2.7cm}p{1cm} p{1.5cm}}
\hline
Equation & Boundary conditions & Unknown & coordinate \\
\hline
\hline
$\frac{d^2y}{dx^2}-\frac{1}{\sqrt{x}}y^{\frac{3}{2}}(x) =0 $ & $y(0)=1, y(\infty)=0$ & $y$ &$x\in [0, \infty]$ \\\\
\hline
$z\{y\frac{d^2y}{dz^2} + \frac{dy}{dz}\} - y \frac{dy}{dz} -2z^2y^3 = 0$ & $y(0)=1, y(\infty)=0$ & $\sqrt{y}$ &$z = \sqrt{x}$\\\\
\hline
$z\frac{d^2y}{dz^2} - \frac{dy}{dz} -4z^2y^{\frac{3}{2}} = 0$ & $y(0)=1, y(\infty)=0$ & $y$ &$z = \sqrt{x}$ \\\\
\hline
\end{tabularx}
\label{formstable}
\end{table}

One special parameter in Thomas-Fermi equation is the first derivative of the unknown function at the region $y'(0)$. This importance is because of some reasons, as the first one, we can refer to the expansion of $y$ about the region, the expansion of $y$ about the region is as follows\cite{Baker1930}:
\begin{equation}\label{baker}
y(x) = 1+\lambda x + \frac{4}{3}x^{\frac{3}{2}} + \frac{2\lambda}{5} x^{\frac{5}{2}} + \frac{1}{3}x^3 + \frac{3\lambda^2}{70} x^{\frac{7}{2}}+ \dots,
\end{equation}
where $\lambda = y'(0)<0$. On the other hand $y'(0)$ can be used to obtain the energy of a neutral atom by the following formula:
\begin{equation}
E = \frac{6}{7}(\frac{4\pi}{3})^{\frac{2}{3}}Z^{\frac{7}{3}}y'(0),
\end{equation}
where $Z$ is the nuclear charge\cite{Laurenzi1990}.

As mentioned above Thomas-Fermi equation has special significance in theoretical physics and thanks to this importance many researchers develop various numerical algorithms to approximate solution of Thomas-Fermi equation. We summarize some previous works in the literature in Table \ref{tab_intro}.\\
\begin{table}
\centering
\caption{A brief bibliography on Thomas-Fermi equation }
\begin{tabularx}{\textwidth}{p{2cm} p{9cm}}
\hline
Years & Description \\
\hline
\hline
1930--1970 & In these years, scientists studied the singularity and convergence of Thomas-Fermi equation, found an analytical solution \cite{Baker1930}, and investigate the asymptotic behavior of $y(x)$ \cite{Sommerfeld1932}.\\
\hline
1970--2000 & Researchers found an alternate analytical solution for Thomas-Fermi equation using perturbative procedure \cite{Laurenzi1990}, solved Thomas-Fermi equation by standard decomposition method \cite{Adomian1998}, Adomian decomposition method and Pad\'{e} approximation \cite{Wazwaz1999,Epele1999}.\\
\hline
2000-2010 & In this decade, scientists proposed various approaches for approximating the solution of Thomas-Fermi equation such as a combination of semi-inverse scheme and the Ritz method \cite{He2003}, piecewise quasilinearization technique \cite{Ramos2004}, an iterative approach and the sweep method \cite{Zaitsev2004}, computing the potential slope at the origin by exploiting integral properties of the Thomas-Fermi equation\cite{Iacono}, rational Chebyshev collocation method \cite{Parand2009}.\\
\hline
2010--2015 & Scientists used semi-analytical and numerical approaches to solve Thomas-Fermi equation with a high accuracy. These techniques are improved Adomian decomposition method \cite{Ebaid2011}, optimal parametric iteration method \cite{Marinca2011}, combination of three schemes based on Taylor series, Pad\'{e} approximates and conformal mappings \cite{Abbasbandy2011}, the Hankel-Pad\'{e} method \cite{Fernandez2011}, an adaptive finite element method based on moving mesh \cite{Zhu2012}, Homotopy analysis method and Pad\'{e} approximates \cite{Turkyilmazoglu2012}, Newton-Kantorovich iteration and collocation approach based on rational Chebyshev functions \cite{Boyd2013}, Sinc-collocation method \cite{Parand2013}, Rational second-kind Chebyshev pseudospectral technique \cite{Kilicman2014}, collocation method on Hermite polynomials \cite{Bayatbabolghani2014}\\
\hline
2015--2018 & Recently, researchers proposed fractional order of rational orthogonal \cite{Parand20171,Parand20172} and non-orthogonal functions \cite{Parand20162,Parand20161} for approximating Thomas-Fermi equation. In 2018, Sabir et al suggest an artificial neural network\cite{Sabir2018} to solve that. Moreover, some other researchers  study coordinate transformations \cite{Zhang2018} for approximating Thomas-Fermi equation and found highly accurate solution to 60 decimal places for $y'(0)$\\
\hline
\end{tabularx}
\label{tab_intro}
\end{table}

\section{Fractional order of rational Gegenbauer (FRG) functions}
\paragraph{}
There are various types of orthogonal polynomials, which have different behaviors and properties. Choosing a good orthogonal function as a basis which behaves as same as the behavior of the exact solution is challenging problem in spectral methods, because we have not the exact solution. But in some problems such as Thomas-Fermi equation, although we have not the exact solution we have some information about the behavior of the solution.

As mentioned in Eq. (\ref{baker}), $y(x)$ can be expanded by a power series of $x^{\frac{1}{2}}$\cite{Parand20172}.
So if we choose fractional functions as a basis it can be fit to the Baker expansion. Moreover, Thomas-Fermi equation is defined in the semi-infinite domain. One choice for the semi-infinite domains is rational functions. Therefore, we select the fractional order of rational Gegenbauer functions as a basis. In this section,  we introduce Gegenbauer polynomials and FRG functions, then we explain how to use FRG function for function approximating.
\subsection{Gegenbauer polynomials}
\paragraph{}
In this paper we use the fractional order of rational Gegenbauer function, where this function is obtained of Gegenbauer polynomial. The Gegenbauer polynomial of degree $n$, $G_n^a(x)$, and order $a>-\frac{1}{2}$  is solution of following differential equation:

\begin{equation}
(1-x^2)\frac{d^2y}{dx^2}-(2a+1)x\frac{dy}{dx}+n(n+2a)y = 0
\end{equation}
where $n$ is a positive integer.

The standard Gegenbauer polynomial $G_n^a(x)$,  is defined as follows:\\
\begin{equation}
G_n^a(x)=\sum_{j=0}^{\lfloor\frac{n}{2}\rfloor}(-1)^j\frac{\Gamma(n+a-j)}{j! (n-2j)! \Gamma(a)}(2x)^{n-2j},
\end{equation}
where $\Gamma(.)$ is the Gamma function.\\
The Gegenbauer polynomials are orthogonal over the interval $[-1,1]$ with the weight function $w(x)=(1-x^2)^{a-\frac{1}{2}}$ which means:
\begin{equation}
\label{orthog}
\int_{-1}^{1}G_n^a(x)G_m^a(x) w(x)dx=\frac{\pi
2^{1-2a}\Gamma(n+2a)}{n!(n+a)(\Gamma(a))^2}\delta_{nm},
\end{equation}
where $\delta_{nm}$ is the Kronecker delta function.

In addition, Gegenbauer polynomials can be obtained by the following recursive formula:
\begin{eqnarray}
G_0^a(x) = 1, ~~~~~ G_1^a(x) = 2ax,~~~~~~~~~~~~~~~~~~~~\\
G_{n+1}^a(x) = \frac{1}{n+1}[2x(n+a)G_n^a(x) - (n+2a-1)G_{n-1}^a(x)], ~~~n\geq1
\end{eqnarray}
\subsection{Fractional order of rational Gegenbauer (FRG) functions}
\paragraph{}
Scientists have been proposing the fractional order of rational functions such as rational Chebyshev \cite{Parand20172}, rational Jacobi \cite{Parand20171}, rational Euler \cite{Parand20161}, etc. to solve some ODEs. A fractional order of the rational Gegenbauer ($FRG$) functions are defined as follows:
\begin{equation}
FRG_n^a(L,\alpha ,x)=G_n^a(\frac{x^{\alpha}-L}{x^{\alpha}+L}),
\end{equation}
in which $L$ and $\alpha$ are real positive numbers. $FRG$ functions are orthogonal functions in semi-infinite interval same as Eq. (\ref{orthog}) according to the weight function $w(x)=(1-(\frac{x^{\alpha}-L}{x^{\alpha}+L}))^{a-\frac{1}{2}}\frac{2\alpha Lx^{\alpha-1}}{(x^{\alpha}+L)^2}$ :
\begin{equation}
\int_{0}^{\infty}FRG_n^a(L,\alpha ,x)FRG_m^a(L,\alpha ,x) w(x)dx=\frac{\pi
2^{1-2a}\Gamma(n+2a)}{n!(n+a)(\Gamma(a))^2}\delta_{nm}.
\end{equation}

\subsection{Approximation of functions}

\begin{definition}
Consider $\Gamma = \{x|0\leq x \leq \infty\}$ and $L_w^2(\Gamma) =\{f:\Gamma\longrightarrow\Re| f~is~measurable~and~||f||_{w}<\infty \} $ where,
\begin{equation*}
  w(x)=(1-(\frac{x^{\alpha}-L}{x^{\alpha}+L}))^{a-\frac{1}{2}}\frac{2\alpha Lx^{\alpha-1}}{(x^{\alpha}+L)^2},
\end{equation*}
and
\begin{equation*}
  ||f(x)||_{w} = \Bigg(\int_{0}^{\infty}f^2(x)w(x)dx\Bigg)^{\frac{1}{2}},
\end{equation*}
is the norm induced by the inner product of the space
\begin{equation*}
  \langle f(x), g(x)\rangle_{w} = \int_{0}^{\infty}f(x)g(x)w(x)dx.
\end{equation*}
\end{definition}

Any function $y(x)\in C(0, \infty)$ can be expanded as the follows:
\begin{equation}
y(x) = \sum_{n=0}^{\infty}a_n FRG_n^a(L,\alpha ,x),
\end{equation}
where
\begin{equation}
a_i = \langle y(x), FRG_i^a(L,\alpha ,x) \rangle = \langle \sum_{n=0}^{\infty}a_n FRG_n^a(L,\alpha ,x), FRG_i^a(L,\alpha ,x) \rangle,
\end{equation}
that is,

\begin{equation}
a_n = \frac{n!(n+a)(\Gamma(a))^2}{\pi 2^{1-2a}\Gamma(n+2a)} \int_0^{\infty}FRG_n^a(L,\alpha ,x) y(x)w(x)dx,
\end{equation}

Now let assume
\begin{equation*}
V_m = span\{FRG_0^a(L,\alpha ,x), FRG_1^a(L,\alpha ,x), \dots, FRG_m^a(L,\alpha ,x)\},
\end{equation*}
is a finite dimensional subspace, therefore $V_m$ is a complete subspace of $L_w^2(\Gamma)$\cite{Boydbook,Fox1968,Parand20172}. Let define the $L_w^2(\Gamma)$-orthogonal projection $\Pi_{N, w} : L_w^2(\Gamma) \rightarrow V_m$, that for any function $y \in L_w^2(\Gamma)$:
\begin{equation}
  \langle \Pi_{N, w}y - y, v\rangle = 0, ~~~ \forall v \in V_m.
\end{equation}

It is clear that $\Pi_{N, w}y$ is  the best approximation of $y(x)$ in $V_m$ and can be expanded as\cite{Goubook}:
\begin{equation}
\Pi_{N, w}y =  y_m(x) = \sum_{i=0}^{m} a_i FRG_i^a(L,\alpha ,x).
\end{equation}


\section{Application of the methods}
\paragraph{}

In this section, two approaches based on Newton method and spectral collocation algorithm are explained to approximate the solution of Thomas-Fermi equation. In one of them, we use the Newton method to linearize the Thomas-Fermi equation and then solve the several linear ODEs by spectral method that we call this method pre-Newton method. In other method we convert the Thomas-Fermi equation to a nonlinear system of algebraic equation by using spectral algorithm, then, solve this nonlinear system by using classical Newton method which we call this method post-Newton. These two schemes are illustrated as follow.


\subsection{Pre-Newton approach for Thomas-Fermi equation}
\paragraph{}
Solving system of nonlinear algebraic equations by using traditional Newton type solvers have three
major practical difficulties. The first one is selecting start point which yields the convergence
of the iterations. The second one is computing the Jacobian matrix of the system of equations at
each iteration that has a lot of computational load to the algorithm. The last one is inverting
a Jacobian matrix at each iteration which is the most expensive step of the algorithm. In the
post-Newton approaches for solving nonlinear ODEs, we should overcome these difficulties\cite{Boydbook}. In order to avoid
these difficulties, we can apply the Newton method directly to the nonlinear ODE. In the next part, a famous
Newton-type algorithm is described which converts nonlinear ODEs to a sequence of linear differential
equations.
\subsubsection{Newton--Kantorovich method }
\paragraph{}
Newton--Kantorovich method is a well-known and strong approach to convert nonlinear ODEs to linear ones which was introduced by
Bellman and Kalaba \cite{Bellman,Conte1981,Ralston1988}. This approach obtains the solution of a nonlinear ODE by solving a
sequence of linear differential equations \cite{Moayeri}. In fact, approximating the solution of a
nonlinear equation is more complicated than a linear one; therefore, by using Newton--Kantorovich method, the solution of
the sequence of the linear differential equations converges to the solution of the original
nonlinear ODE\cite{Bellman,Mandelzweig1999,Mandelzweig2001}. This method is based on approximating a nonlinear function by using linear part of Taylor expansion of that function.

%

This fact can be extended to linearize a nonlinear ODE. In order to show how Newton--Kantorovich method works we consider a $n$-th order nonlinear ODE over the interval $[0,b]$ as follows\cite{Mandelzweig2001}:

\begin{equation}\label{generaleq}
L^{(n)} y(x) = f(y(x), y^{(1)}(x), \dots, y^{(n-1)}(x), x),
\end{equation}
with the following boundary conditions:
\begin{equation}
B_k(y(0), y^{(1)}(0), \dots, y^{(n-1)}(0)) = 0 ~~~~ k = 1, 2, \dots, l,
\end{equation}
and
\begin{equation}
B_k(y(b), y^{(1)}(b), \dots, y^{(n-1)}(b)) = 0 ~~~~ k = l+1, l+2, \dots, n,
\end{equation}
where $L^{(n)}$ is a linear $n$-th order ordinary differential operator and $f$ and $B_1, B_2, \dots, B_n$ are nonlinear functions of $y(x)$ and its $n-1$ derivatives $y^{(s)}, s= 1, 2, \dots, n-1$. If we apply Newton--Kantorovich method on Eq. (\ref{generaleq}) the $(r+1)$-th iterative approximation of $y(x)$ is obtained by solving follow linear ODE,
\begin{equation}
\begin{split}
    L^{(n)} y_{r+1}(x) = f(y_{r}(x), y_{r}^{(1)}(x), \dots, y_{r}^{(n-1)}(x), x) +~~~~~~~~~~\\ \sum_{s = 0}^n \big(y_{r+1}^{(s)}(x) - y_{r}^{(s)}(x)\big) f_{y^{(s)}}(y_{r}(x), y_{r}^{(1)}(x), \dots, y_{r}^{(n-1)}(x), x),
\end{split}
\end{equation}
where $y^{0}_r(x)$ is a notation for $y_r(x)$. Also, the linearized boundary conditions are obtained as follows:
\begin{equation}
\sum_{s = 0}^{n-1} \big(y_{r+1}^{(s)}(0) - y_{r}^{(s)}(0)\big) B_{k y^{(s)}}(y_{r}(0), y_{r}^{(1)}(0), \dots, y_{r}^{(n-1)}(0), 0) = 0, ~~~~~~ k = 1, \dots, l,
\end{equation}
and
\begin{equation}
\sum_{s = 0}^{n-1} \big(y_{r+1}^{(s)}(b) - y_{r}^{(s)}(b)\big) B_{k y^{(s)}}(y_{r}(b), y_{r}^{(1)}(b), \dots, y_{r}^{(n-1)}(v), b) = 0, ~~~~~~ k = 1, \dots, l.
\end{equation}
It is worth to mention that in the above formulas $f_{y^{(s)}} = \frac{\partial f}{\partial y^{(s)}}$ and $B_{k y^{(s)}} = \frac{\partial B_k}{\partial y^{(s)}}$ for $s = 0, 1, \dots, n-1$

By implementing Newton--Kantorovich method on Eq. (\ref{TFE}), the $(i+1)-th$ iteration linear ODE for approximating the
solution of Thomas-Fermi equation is as follows ($i = 0, 1, 2,\dots$):
\begin{equation}
\label{TFE-NKM}
\sqrt{x}y''_{i+1}(x)-\frac{3}{2}\big(y_{i}(x)\big)^{\frac{1}{2}}y_{i+1}(x)=\frac{-1}{2}\big(y_i(x)\big)^{\frac{3}{2}
},
\end{equation}
with the following boundary conditions:
\begin{eqnarray}
\label{BTFE-NKM}
y_{i+1}(0)=1,~~~~\lim_{x\to\infty} y_{i+1}(x)=0.
\end{eqnarray}
An initial guess $y_0(x)$ is required for the first step of the Newton--Kantorovich method. It is proved that when the initial
guess satisfies one of the boundary conditions, the Newton--Kantorovich method will be convergent\cite{Bellman}.
Thus, we consider $y_0(x)$=1.

\subsubsection{Collocation method in the pre-Newton method}
\paragraph{}
\label{col-QLM}
The spectral collocation method based on FRG functions is applied to Eq. (\ref{TFE-NKM}) at each
iteration. According to the boundary conditions in Eq. (\ref{BTFE-NKM}), we approximate $y_{i+1}(x)$ in
$(i+1)-th$ iteration as:
\begin{equation}
\label{approx}
y_{i+1}(x)\simeq y^N_{i+1}(x)=1 +x \sum_{j=0}^{N-1} a^{i+1}_j FRG_j^a(L,\alpha,x).
\end{equation}
where $a^{i+1}_j$ is the $j-th$ unknown coefficient in $(i+1)-th$ iteration. Equation
(\ref{approx}) satisfies the boundary condition $y(0)=1$. To satisfy the other boundary condition,
we choose a sufficiently large number $K$ and consider $y_{i+1}(K)=0$. The Eq. (\ref{approx}) is
replaced in Eq. (\ref{TFE-NKM}); afterwards, the residual function is obtained:\\
\begin{equation}
\label{res}
Res_{i+1}(x)=\sqrt{x}y''^N_{i+1}(x)-\frac{3}{2}\big(y^N_{i}(x)\big)^{\frac{1}{2}}y^N_{i+1}(x)+\frac{1}{2}\big(y^
N_{i}(x)\big)^{\frac{3}{2}},
\end{equation}
The roots of $FRG^a_N(L,\alpha,x)$ are considered as the collocation points which are collocated in
Eq. (\ref{res}) and a system of linear algebraic equations is established. By solving this system at
each iteration, $y(x)$ is approximated.
\begin{equation}
Res_{i+1}(x_j)=0,~~~j=0,...,N-1.
\end{equation}

\subsection{Post-Newton approach for Thomas-Fermi equation}
In the post-Newton approach for solving Thomas-Fermi equation, we use a fully spectral technique same collocation method to solve nonlinear equation Eq. (\ref{TFE}) without any linearization method. In this method, by using spectral collocation method based on fractional order of rational Gegenbauer functions, we convert nonlinear Thomas-Fermi equation Eq. (\ref{TFE}) to a system of nonlinear algebraic equations. In this method, the unknown solution $y(x)$ of Thomas-Fermi is approximated by the following series:
\begin{equation}
\label{approx2}
y(x)\simeq y_N(x)=1 +x \sum_{j=0}^{N-1} a_j FRG_j^a(L,\alpha,x).
\end{equation}
Then by substitution $y^N(x)$ instead of  $y(x)$ in  Eq. (\ref{TFE}) the residual function is constructed as follows:
\begin{equation}
\label{res2}
Res(x)=\sqrt{x}y''^{N}(x)-\frac{1}{\sqrt{x}}\big(y^{N}(x)\big)^{\frac{3}{2}}.
\end{equation}
Now, there are $N$ unknown coefficients $a_i,i=0,1,...,N-1$, to find these unknowns we need $N$ equations. By using collocation technique and roots of fractional order of rational Gegenbauer function of order $N$, and by substitution these nodes in residual function we construct $N$ nonlinear equations as follows:
$$F_i=Res(x_i)=0,~~i=0,1,...,N-2,$$
to satisfy the boundary condition in infinite we set $F_{N-1}=y_N(L)=0$ for sufficient large $L$.\par
So, $F:\mathbb{R}^N\rightarrow \mathbb{R}^N$ is a nonlinear function and therefore, finding the solution of Eq. (\ref{TFE}) has been transformed to find the solution of the nonlinear system of equations:
$$F(A)=0, ~~~~A=[a_0,a_1,...,a_{N-1}]^T,$$
Now, we have transformed solving the nonlinear differential equation to finding the root of a nonlinear $\mathbb{R}^N\rightarrow \mathbb{R}^N$ function.
\par
One of the best methods to solve a nonlinear system is the classical Newton iterative method, that by using Tylor expansion:
\begin{eqnarray}
F(x_{n+1})=F(x_n)+(x_{n}-x_{n+1})F'(x_n),
\end{eqnarray}
presuppose $x_{n+1}$ be root of $F(x)$, $F(x_{n+1})=0$ :
\begin{eqnarray}
&& F(x_n)+(x_{n+1}-x_{n})F'(x_n)=0,\\
&& \Rightarrow x_{n+1}=x_n- F'(x_n)^{-1} F(x_n),
\end{eqnarray}
$F'(x)=J(x)$ is the $n\times n$ Jacobian matrix and is defined as follows:
\begin{equation}
J_{ij}=\bigg(\frac{\partial f_i}{\partial x_j}\bigg),
\end{equation}
therefore:
\begin{equation}
x_{n+1}=x_n-J(x_n)^{-1}{F(x_n)}.
\end{equation}
In fact in each iteration, a linear system must be solved:
\begin{equation}\label{NewtonII}
\left\{
\begin{array}{l}
x_{n+1}=x_n+\delta x_n\\\\
J(x_n) \delta x_n=F(x_n).
\end{array}
\right.
\end{equation}
In this paper, we use $LU$ 
method to solve linear system  $J(x_n) \delta x_n=F(x_n)$ in each iteration of Newton method. Initial guess of the post-Newton method is the simple vector $x_0 = (1, 1, \dots, 1)^{T}$.

\section{Numerical results and discussion}
\begin{figure}
\centering
\begin{subfigure}{0.4\textwidth}
\centering
\includegraphics[height=1.8in,keepaspectratio=true]{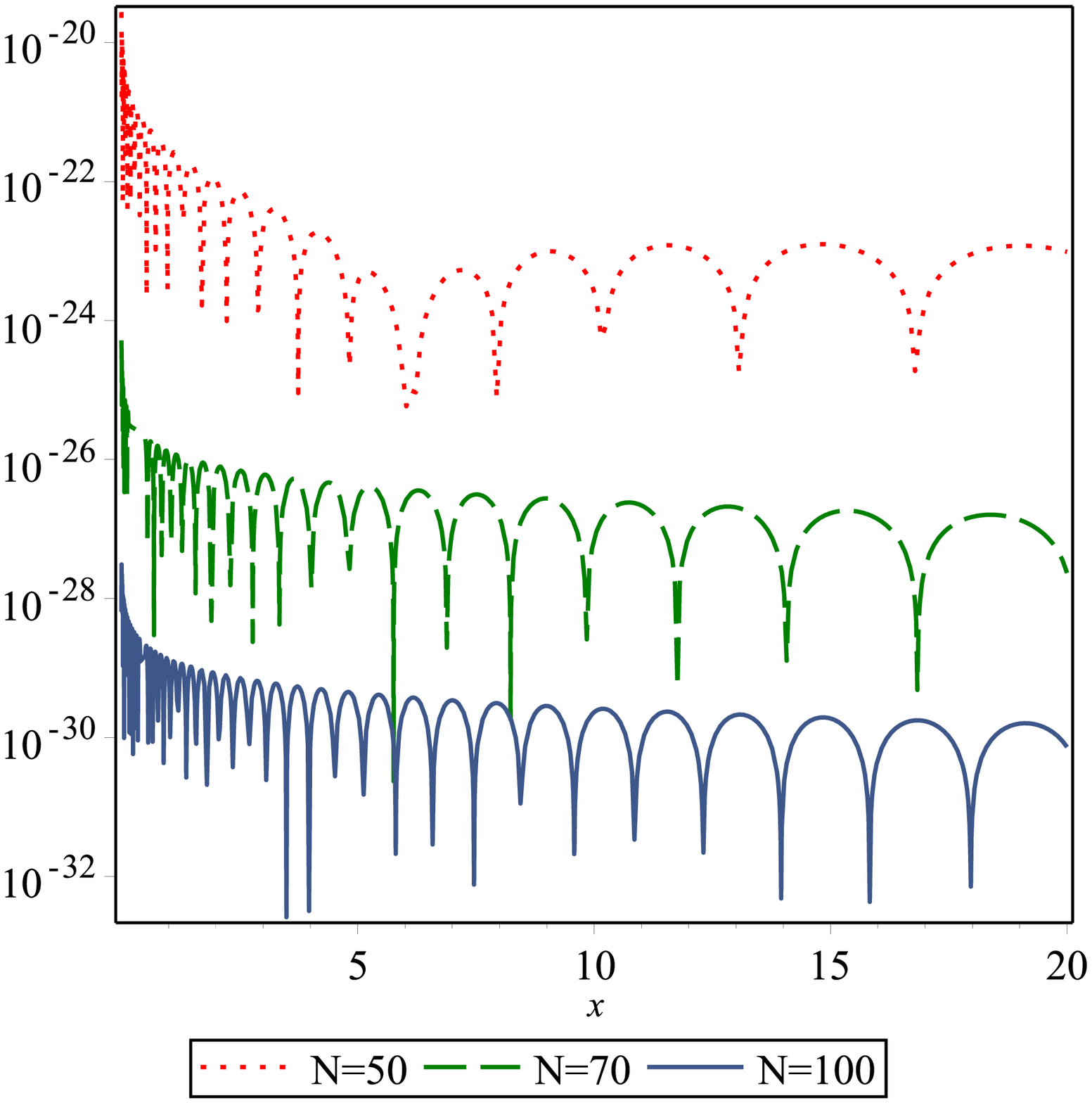}
\caption{}
\end{subfigure}
\begin{subfigure}{0.4\textwidth}
\centering
\includegraphics[height=1.8in,keepaspectratio=true]{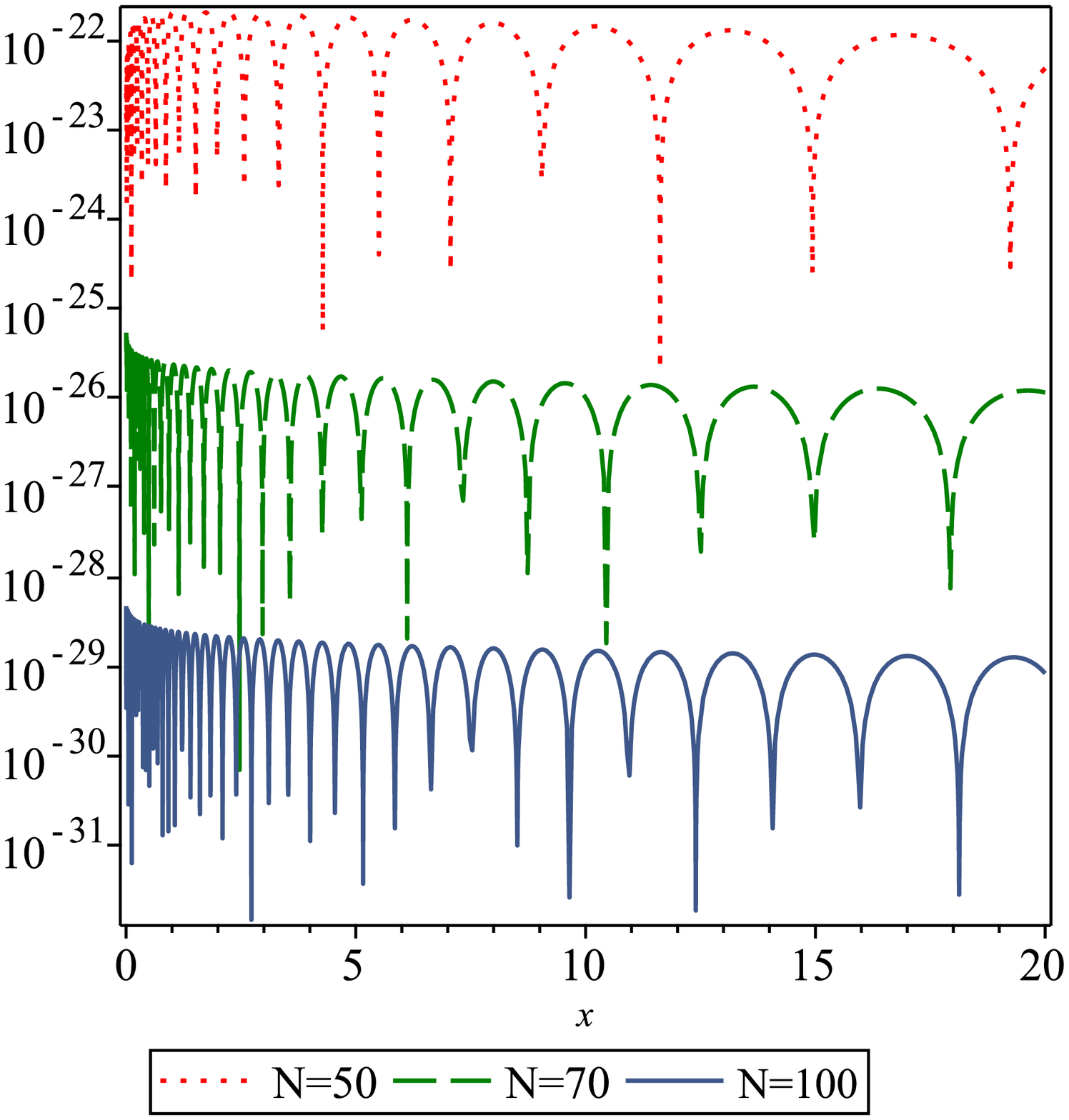}
\caption{}
\end{subfigure}
\caption{Graph of the logarithm of the absolute residual errors for different $N$ for (a) pre-Newton (b) post-Newton.}
\label{ress}
\end{figure}

\paragraph{}
According to Boyd's book \cite{Boydbook}, $L$ can be chosen by "The experimental trial-and-error
method"; so, we consider $L=3$ in pre-Newton and $L = 2.828$ in post-Newton, also we consider
$\alpha=\frac{1}{2},a=\frac{1}{2}$ in the both and report the results. It is worth to mention that
all the computations are done by Maple, in a personal computer with the following hardware
configuration: desktop 64-bit Intel Core i5 CPU, 8GB of RAM, 64-bit Operating System. In
\cite{Zhang2018}, Zhang and Boyd calculated an the approximate solution for  $y'(0)$ with high
accuracy; thus, the results of this study is compared with \cite{Zhang2018} and found that the
obtained results are as accurate as \cite{Zhang2018}. The logarithm of absolute residual error for
Thomas-Fermi equation in the best iteration is represented in Fig. \ref{ress}.This figure shows when the number of
collocation points increases, the residual error tends to the zero. The value of $y'(0)$ is
presented in Table \ref{y'(0)} and compared with the obtained solution by state-of-the-art
methods. Table \ref{grap1} contains the values of $y(x)$ and $y'(x)$ for different values of $x$.\\

\begin{figure}
\centering
\includegraphics[height=2in,keepaspectratio=true]{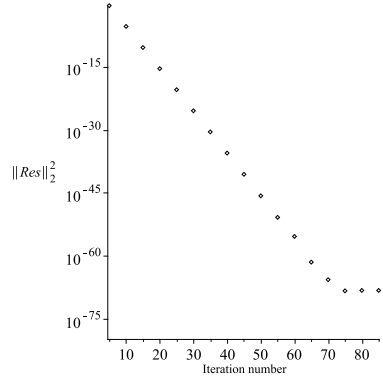}
\caption{Graph of the $\log(||Res||^2_2)$ for the post-Newton method with $N=200$ at different iterations.}
\label{res_final}
\end{figure}

\begin{table}[H]
\small
\centering
\caption{Comparison of the obtained values of $y'(0)$ by some researchers}
\label{y'(0)}
\begin{tabular}{lll}
\hline
Author/Authors & Year & Obtained value of $y'(0)$ \\ \hline
Boyd   \cite{Boyd2013}&(2013)& -1.5880710226113753127186845  \\
Parand et al  \cite{Parand20171}&(2017)& -1.588071022611375312718684509423950109  \\
Parand and Delkhosh (N=300) \cite{Parand20172} &(2017)& -1.58807102261137531271868450942395010951 \\
Zhang and Boyd (N=600)\cite{Zhang2018}&(2018) & -1.588071022611375312718684509423950109452746621674825616765677\\
pre-Newton (N=100) && -1.58807102261137531271868450942395010945274662\\
post-Newton (N=200)  && -1.588071022611375312718684509423950109452746621674825616765677 \\
\hline
\end{tabular}
\end{table}
\begin{table}[H]
\centering
\small
\caption{Values of $y(x)$ and $y'(x)$ obtained by the presented methods for the various values of $x$}
\begin{tabular}{|c|c|c|c|c|c|}
\hline
$y(x)$ and $y'(x)$& $x$  & pre-Newton ($N$=100 and iteration=40) & post-Newton ($N$=200 and iteration=85) \\
\hline
\multirow{4}{*}{$y(x)$}
& 0.5     & 0.6069863833559799094944460701740221017049 & 0.6069863833559799094944460701740842378463\\
& 3     & 0.1566326732164958413398134404775366125433 & 0.1566326732164958413398134404779118302783\\
& 10     & 0.0243142929886808641901103881732913695553 & 0.0243142929886808641901103881763049683685\\
& 50     & 0.0006322547829849047267797787287302055560 & 0.0006322547829849047267797787427886658114 \\
& 200   & 0.0000145018034969457646803986629623432665 & 0.0000145018034969457646803987687276929118 \\
& 5000 &  0.0000000011309267063430848076021125559361 & 0.0000000011309267063430848263855178787850 \\
\hlineB{4}
\multirow{4}{*}{$y'(x)$}
& 0.5     & -0.4894116125745380886470058475611743123609 & -0.4894116125745380886470058475573462887337\\
& 3  & -0.0624571308541209762287048999941581989893 & -0.0624571308541209762287048999995217973789 \\
& 10     &-0.0046028818712692545025435118554873081322 & -0.0046028818712692545025435118515886154232 \\
& 50     &-0.0000324989020482588146242006692476761611 & -0.0000324989020482588146242006802396097650 \\
& 200   &-0.0000002057532316475268926057043855114949  & -0.0000002057532316475268926056858363001742 \\
& 5000 &-0.0000000000006753397121638834659796119395 & -0.0000000000006753397121638835144503744957 \\
\hline
\end{tabular}
\label{grap1}
\end{table}
One of the advantages of the post-Newton approach is its computational speed. This approach is much faster than pre-Newton; as the iterations can be increase to 85 with an acceptable runtime. In Table \ref{tab_run}, pre-Newton and post-Newton methods are compared in runtime with the different number of collocation points and iterations. It is derived that post-Newton is much faster than the other approach; therefore, we can consider a larger number of iterations for the post-Newton than pre-Newton. The logarithm of $||Res||^2$ at different iterations of the post-Newton method for Thomas-Fermi equation by using $200$ points is represented in Fig. \ref{res_final}.

\begin{table}
\centering
\small
\caption{Runtime for the proposed methods with the different values of $N$ and iteration}
\begin{tabular}{|c|c|c|c|}
\hline
$N$ & Iteration & Runtime for pre-Newton (s) &  Runtime for post-Newton (s) \\
\hline
\multirow{3}{*}{50}
& 20     & 39.901 & 32.392 \\
& 30     & 60.626 & 48.688 \\
& 40   & 80.361 & 56.359\\
\hlineB{3}
\multirow{3}{*}{70}
& 20     &104.191 & 76.690 \\
& 30     &172.292 & 111.274 \\
& 40   &211.721 & 136.428 \\
\hlineB{3}
\multirow{3}{*}{100}
& 20     &343.672 & 208.404 \\
& 30     &469.262 &304.123 \\
& 40   &622.255  &392.937 \\
\hline
\end{tabular}
\label{tab_run}
\end{table}

\section{Conclusion}
\paragraph{}
In this paper, we introduced and compared two point of views to solve nonlinear boundary problems over the semi-infinite interval. These two approaches are called pre-Newton method and post-Newton method, respectively.
The pre-Newton method is based on applying Newton--Kantorovich algorithm to the nonlinear ODE and solving the obtained linear ODEs from Newton--Kantorovich method by using collocation algorithm.
The post-Newton method is based on applying collocation algorithm directly to the nonlinear ODE and then solve the obtained nonlinear system of algebraic equations by classical iterative Newton method. The collocation algorithm which is used is based on orthogonal functions in the interval $[0,\infty)$ which are called the fractional order of the rational Gegenbauer. Since the significance of the Thomas-Fermi equation, here, we consider it as a test problem. In the Thomas-Fermi equation the value of $y'(0)$ has important information in physics and scientists attempt to approximate that precisely. Therefore, we compare the approximation solution in $y'(0)$ with the other numerical methods and realize that our proposed method is effective. The approximate solutions for $y(x)$ and $y'(x)$ for various values of $x$ are represented. Additionally, the suggested methods are compared in runtime to find out which method is more efficient. According to the results, the post-Newton approach is faster and more accurate than the pre-Newton approach. It is worth to mention that one of the limitations of the proposed algorithms is ill-posedness of systems of algebraic equations. This limitation causes that we can not increase the number of collocation points.

 %

\end{document}